# A global optimization approach for antenna design using analytical derivatives from high-frequency simulations

Thomas Most, Ansys Dynardo GmbH, Weimar, Germany

Peter Krenz, Ralf Lampert, Ansys Inc., Canonsburg, PA, United States

## 1  Summary

Antennas are more prevalent than ever enabling 5G connectivity for wide ranging applications like cellular communication, IoT, autonomous vehicles, etc. Optimizing an antenna design can be challenging and employing traditional optimization techniques oftentimes require evaluating ~100s of design variations, which is time prohibitive to tight design deadlines. Instead, the common optimization approach relies on engineering expertise and solving only a handful design variations to manually determine a design with good-enough performance.

Here, we present a novel approach for antenna optimization, which utilizes analytical derivatives calculated by Ansys HFSS to dramatically reduce the number of necessary simulation runs. In the presented approach a local updating scheme, which converges quickly to the next local optimum, is extended to a global search strategy intended to detect several local and global optima. This global search procedure uses a consistent interpolation model to approximate the simulation model outputs and its derivatives depending on the design parameters values. The presented approach is implemented in an optimization workflow, which automatically interacts with the HFSS simulation model. Since this optimizer relies on the analytical derivatives calculated by HFSS, it is applicable to any design optimization where geometry or material properties are varied, and the device response is described by SYZ-parameters or far-field quantities. The usability of this new optimizer is therefore wide ranging from antenna design to signal integrity applications.

## 2  Antenna design using Ansys HFSS

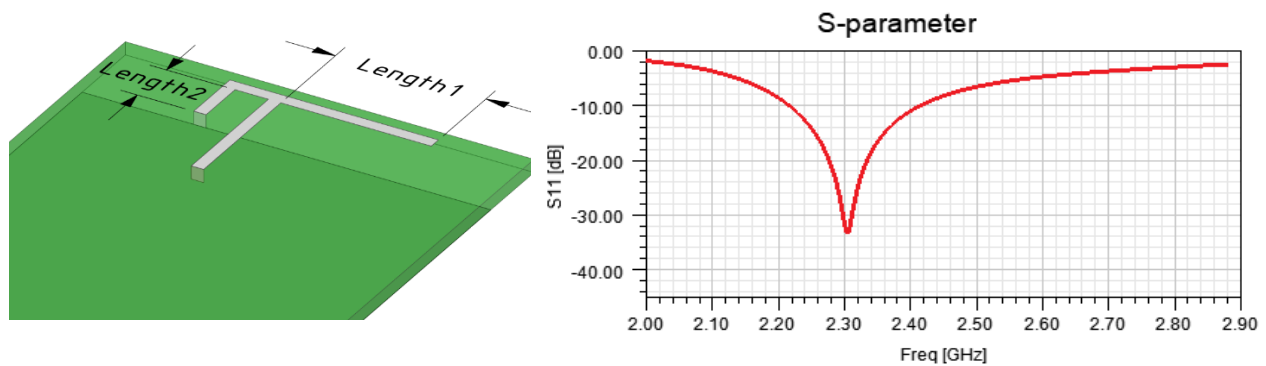

Fig. 1: Single band antenna with calculated return loss frequency spectrum

Ansys High-Frequency Structure Simulator (HFSS) is a general-purpose full-wave 3D electromagnetic (EM) simulation software for simulating and optimizing high-frequency electronics products like antennas, antenna arrays, high-speed interconnects, and printed circuit boards to name a few. Using Ansys HFSS allows engineers to accurately evaluate the performance of complex designs before the prototype phase. Many antenna design application calculations cannot be done by hand, making this high-frequency software extremely valuable for the end-user.

Ansys HFSS can be used for antenna design to calculate the return loss spectrum based on a given geometry. In Figure 1 a single-band antenna with resulting magnitude of the return loss spectrum is shown. Additionally, HFSS has the option to compute analytical derivatives of the return loss with respect to input parameters that control geometrical aspects or material properties of the antenna model [6]. In this paper, we will present an optimization approach, which utilizes this additional solver information very efficiently to find an optimal antenna design.

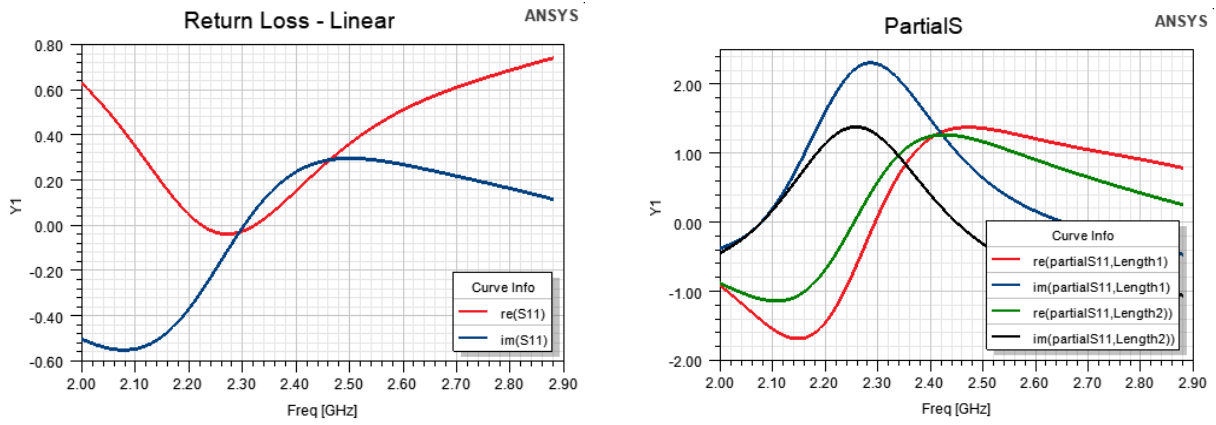

Fig. 2: Real and imaginary parts of the return loss of the single band antenna and its derivatives with respect to the geometry parameters

## 3     Local optimization approach

In a first step, we use the analytical derivatives of the real and imaginary signals responses for a local approximation model based on a linear Taylor series approach:

$$S_{real,local} = S_{real,nominal} + \sum \frac{\partial S_{real,nominal}}{\partial \text{input}_i} \Delta \text{input}_i$$

$$S_{imag,local} = S_{imag,nominal} + \sum \frac{\partial S_{imag,nominal}}{\partial \text{input}_i} \Delta \text{input}_i$$

Based on that assumption, the real and imaginary signal parts are approximated for each discretization point linearly using a single evaluation point. From these approximations, the signal values on the dB scale can be assembled as:

$$|S_{local}|(dB) = 20 \log_{10} \left| \sqrt{\left(S_{real,local}\right)^2 + \left(S_{imag,local}\right)^2} \right|$$

As shown in Figure 3, the linear approximations of the real and imaginary parts may result in significant non-linear functions for the approximated values on the dB scale.

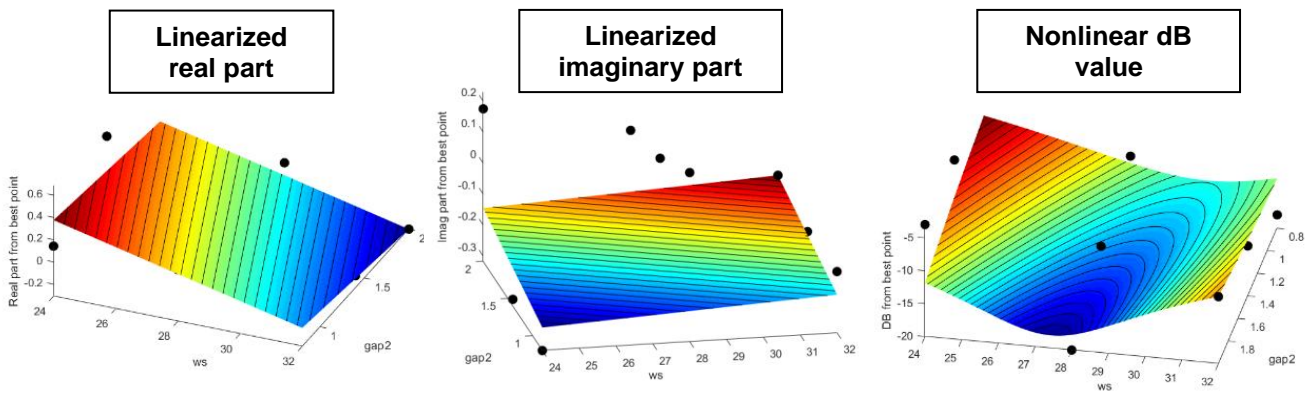

Fig. 3: Linearized real and imaginary parts and assembled nonlinear dB value

Starting at an initial design point, the local approximation is evaluated in a defined region around that evaluation point and a new candidate point is obtained with a local optimization search using only the approximated objective function instead of real solver runs, where the signal responses of the approximated signals at the dB scale are considered. For that new optimal candidate point, the HFSS simulation is executed again and the nominal signal values and their derivatives are updated. If this new design shows a better performance than the previous evaluation point, the search domain is adapted for the next iteration around the new point and reduced in size by a given factor as shown in Figure 4. If the start point is selected well, this iteration converges usually within 5-10 steps to a significantly improved design. In Figure 4, the convergence is shown for the single-band antenna where the objective was to minimize the return loss at 2.44 GHz.

The local approach works well, if only one dominant local optimum exists in the vicinity of the initial design. If the objective contains several local minima, the optimization procedure may get stuck in one of these depending on the initial start design. This problem is shown in Figure 4 for a well-known dual band antenna example [Gai2010], whereas for illustration purposes only the two parameters *ws* and *gap2* are varied. As optimization objective, the minimization of the return loss at 5.6 GHz was considered. In Figure 4 it is illustrated, how the local optimization approach converges in the 2D design space. Depending on the start design, some of the optimization runs did not find the global best solution.

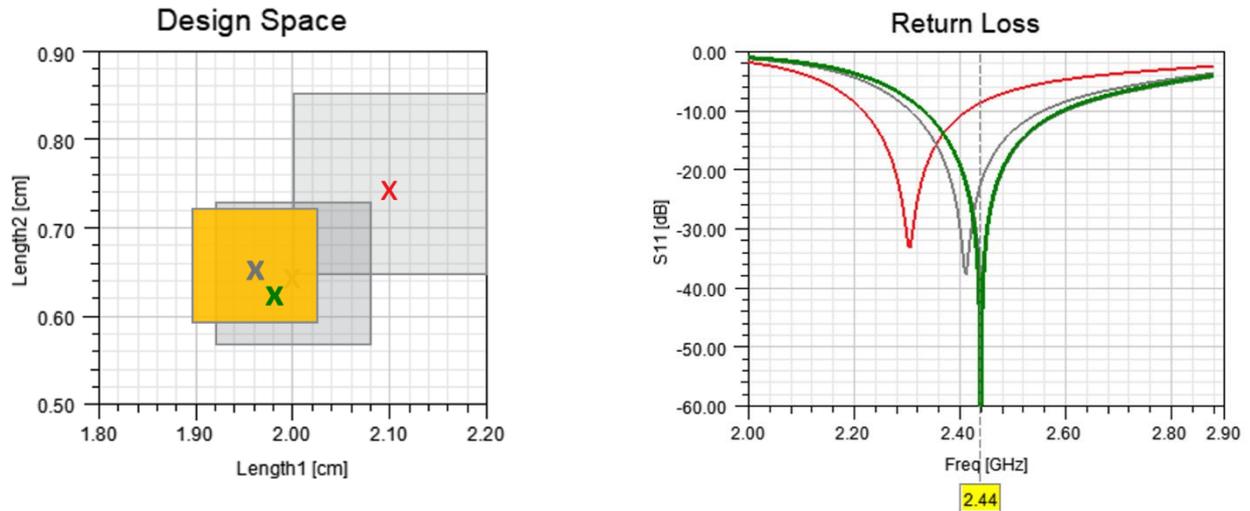

Fig. 4: Local optimization approach using linear approximations of the real and imaginary return loss signals within a shrinking search domain

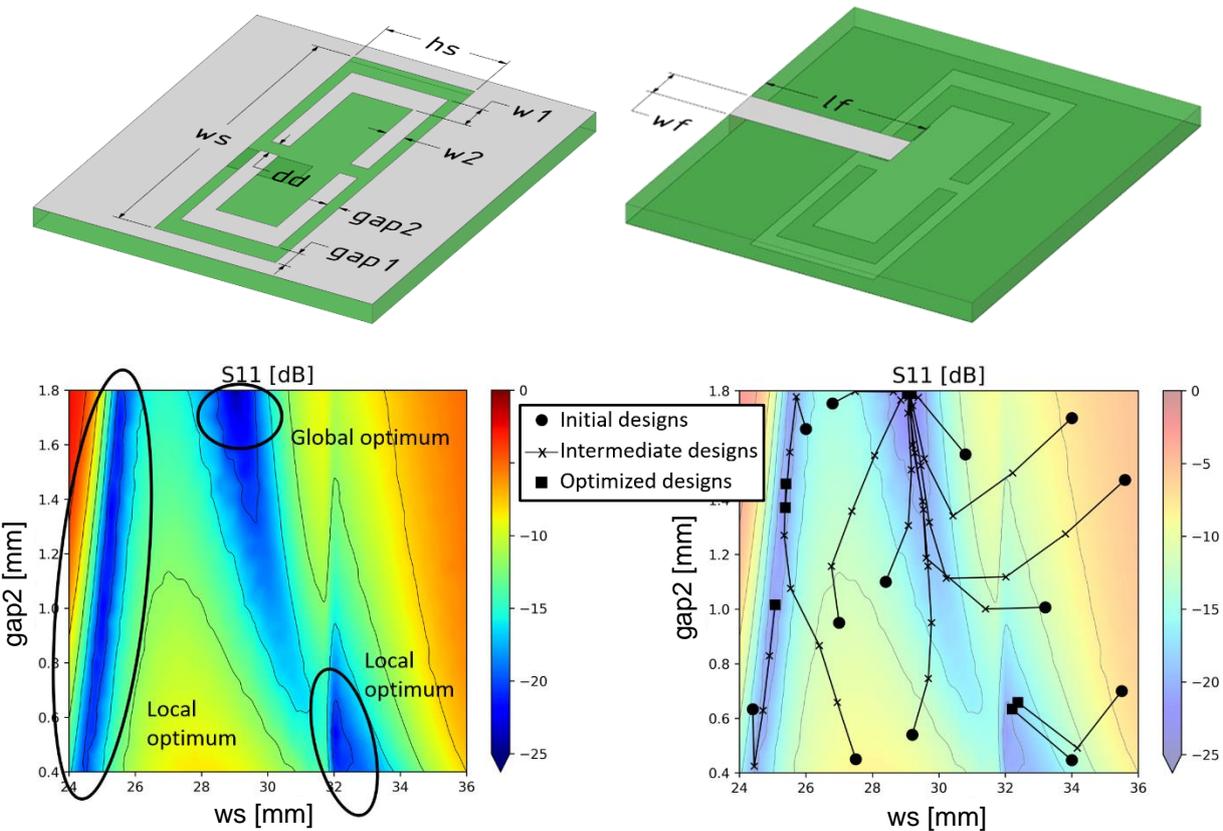

Fig. 5: Local and global optima and convergence of several local optimization runs for the dual-band antenna example optimized for a single resonance using only *ws* and *gap2* as design parameters

## 4    Global optimization approach

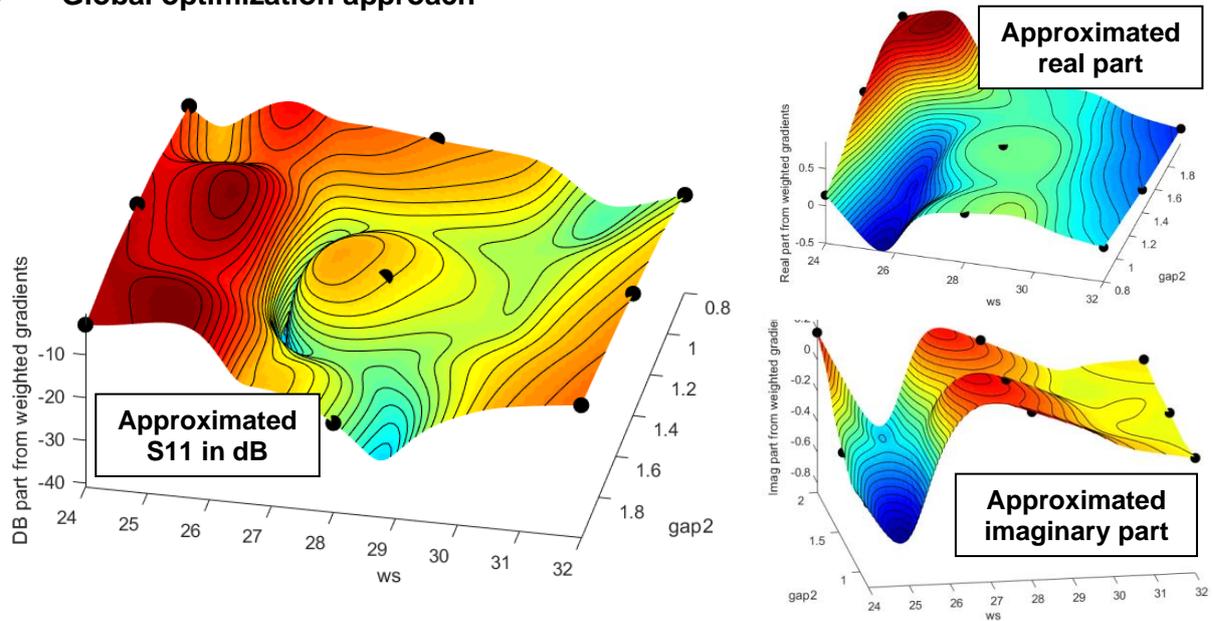

Fig. 6: Approximated real and imaginary parts and their assembly in a global approximation for the S11 response of the 2D dual band antenna in the dB scale

In order to overcome the limitation of the local approach, a novel global optimization approach is introduced, where the partial derivatives of the design points are used to build a global approximation model. Based on an initial design set, which is typically an efficient Design of Experiment scheme [Myers2002] such as a full factorial design or Latin Hypercube Sampling [Hungtinton1998], a global approximation model is generated using a weighted interpolation of the individual linearized approximations of all available design points:

$$S_{real,global} = \sum w_{point\ i} \cdot S_{real,local,point\ i}; \quad S_{imag,global} = \sum w_{point\ i} \cdot S_{imag,local,point\ i}$$

As weighting function $w$, the interpolating Moving Least Squares approximation [Most2005] is utilized, since it represents the design point values exactly. As a result, the approximated responses for the real and imaginary parts and all derived outputs agree exactly with the simulation results at each design point used for the approximation. Similar to the linearized approach, the real and imaginary parts of all investigated response signals are approximated first and then the dB scale and the derived objective values are calculated from these approximations. In Figure 6 the approximation functions are shown for the dual-band example based on an initial 3x3 point full factorial design. It can be seen in the figure, that the approximated real and imaginary parts are quite smooth functions, whereas the assembled dB values show one significant resonance domain.

Based on this global approximation, a refinement-strategy is introduced which efficiently generates new designs in regions, where the global optimum is expected. For this purpose, the Expected Improvement strategy according to [Jones] is applied:

$$E[I(obj)] = \left(obj_{best} - obj_{approx}\right) \cdot \Phi\left(\frac{obj_{best} - obj_{approx}}{\sigma_{obj}}\right) + \sigma_{obj} \phi\left(\frac{obj_{best} - obj_{approx}}{\sigma_{obj}}\right)$$

This approach estimates the possible expected improvement of an approximated single-objective optimization goal function based on a global approximation in the design space, whereas $obj_{best}$ is the best objective value obtained in the previous iterations and $obj_{approx}$ is the approximated objective value. Assuming that the approximated objective is normally distributed, the cumulative density function $\Phi$ and its derivative $\phi$ can be evaluated, if the mean and its standard deviation $\sigma_{obj}$ are known. The mean is given here directly as the approximated scalar objective value based on the real and imaginary part approximations. However, an estimate of the standard deviation of an arbitrary formulation of the objective seems not possible in closed form. Therefore, an estimator based on the difference between the proposed global and the local approximation using the closest available data point in the design space is introduced, which is shown in Figure 7. The candidate point with the maximum expected improvement is chosen for the update of the designs in the next iteration. Since the global and the local approximations agree exactly with the original response values at each considered design point, the estimated error and thus the expected improvement is zero and no new design points are generated too close to the existing data points.

Simultaneous to the global design point update in each iteration, the local approach based on the best design in the iteration history is applied to improve the local convergence. As convergence criterion, a stagnation iteration number is used. In Figure 8 the approximation function of the initial and the final iteration after 15 global and local adaptation steps are shown. The figure indicates, that all local optima and the global optimum are determined by the presented approach.

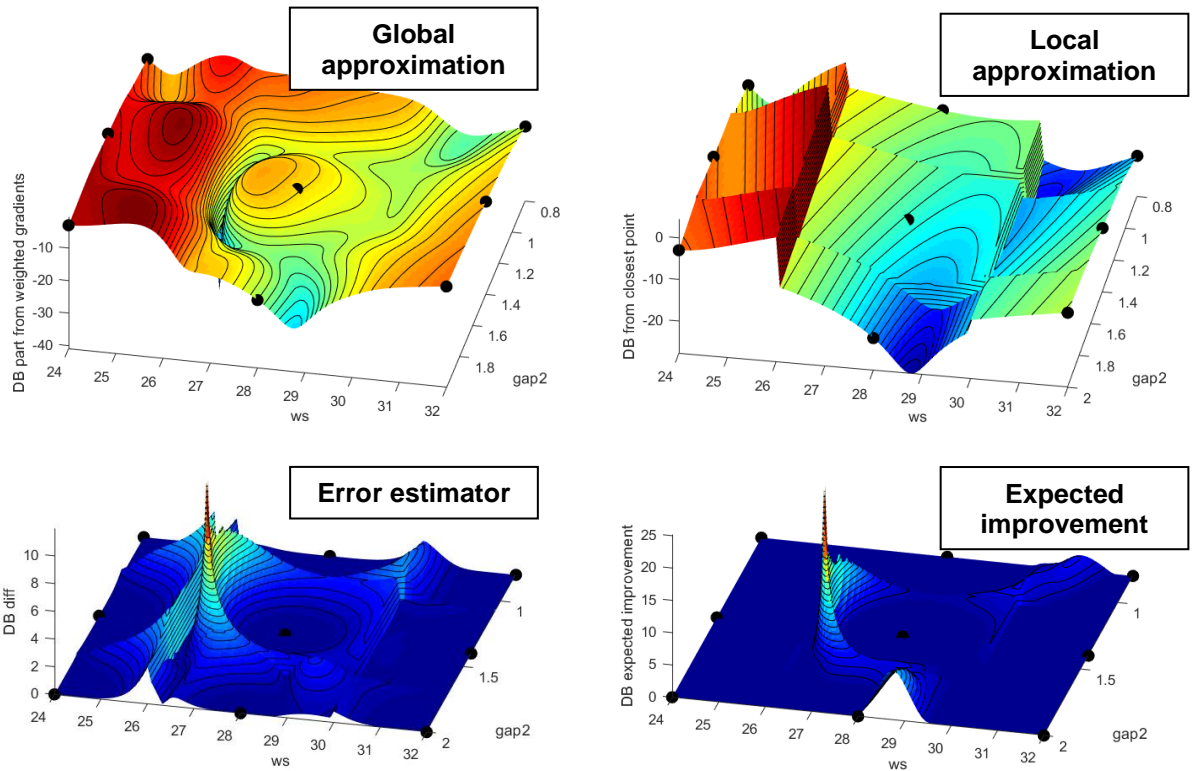

Fig. 7: Approximated objective function and corresponding error estimator for the 2D dual band antenna in the dB scale for the initial design set

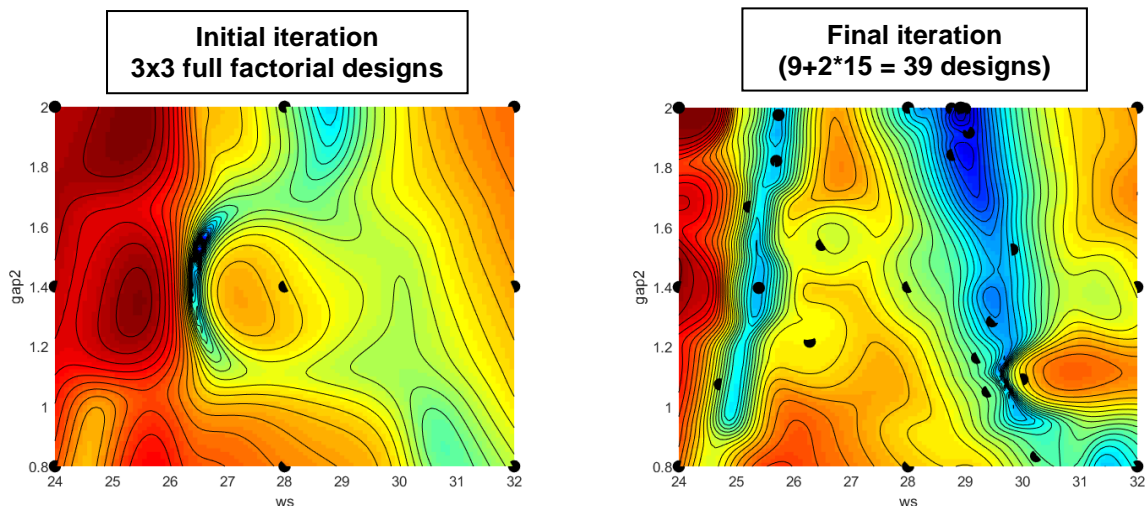

Fig. 8: Initial approximation of the objective function for the 2D dual band antenna and its update with the expected improvement and local refinement approach within 15 iterations

## 5 Dual-band antenna optimization

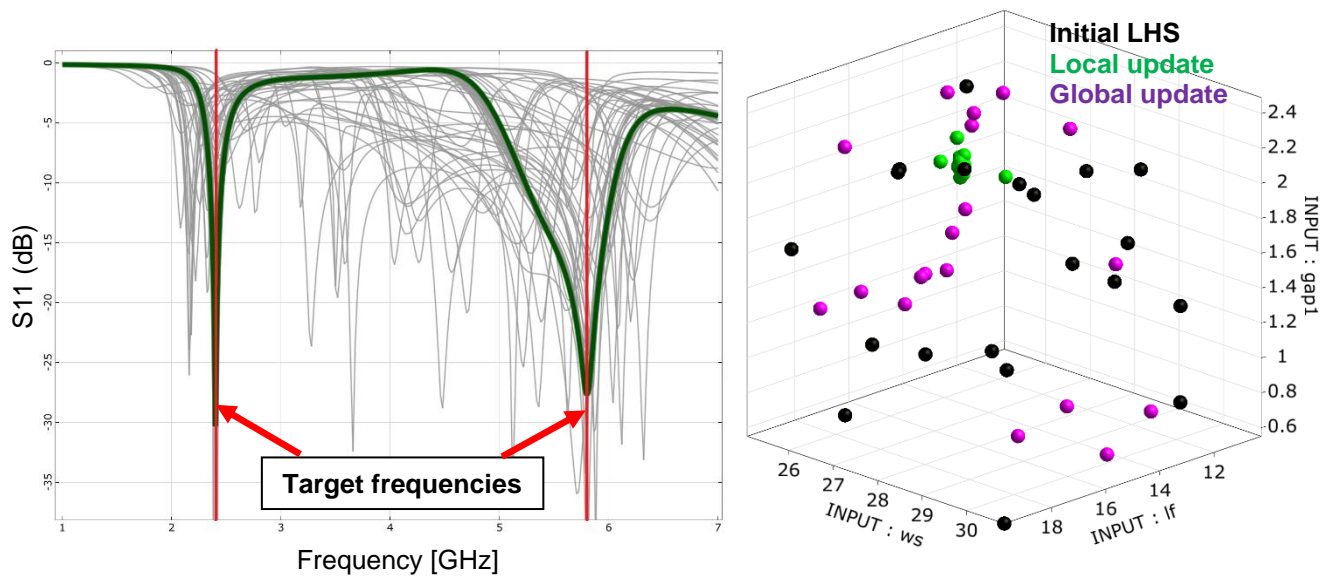

Fig. 9: Dual-band antenna optimization considering 2 resonances and 9 design parameters: designs from local and global refinement procedures and optimal design frequency spectrum

Within this example, the original dual-band antenna example [Gai] introduced in section 3 is optimized regarding the two target WiFi frequencies 2.4 and 5.8 GHz. All 9 geometry parameters shown in Figure 5 are considered as optimization parameters. The maximum dB value from both target frequency values is taken as the objective function which has to be minimized during the optimization procedure. Starting with an initial Latin Hypercube Sampling of 20 designs covering the whole design space, in each iteration the local and global approach is applied to determine two new design candidates. After 20 iterations with total of 40 designs the optimization procedure converged to a very good solution as shown in Figure 9. Additionally, the distribution of the refinement designs is shown in the subspace of the 3 most important parameters. This figure shows, that only one dominant region is detected and thus just one local optimum exists for that example.

## 6 Conclusions

In this paper, we presented a novel optimization approach for antenna design using Ansys HFSS simulation capabilities. Based on a local linearization technique using analytical derivatives from the simulation, a global approximation approach with an efficient refinement was derived. The resulting approximation model is able to detect several local and global optima in the objective function with a very limited number of simulation runs. The presented algorithm is implemented in Ansys optiSLang and is available in the 2022 R2 release. An application of optiSLang and HFSS for a diplexer safety analysis can be found in [7].